\documentclass{amsart}

\usepackage[utf8x]{inputenc}

\usepackage{subfigure}
\usepackage[english]{babel}
\usepackage{amsmath, amsfonts, amssymb}
\usepackage{color}
\usepackage{tikz}
\usetikzlibrary{positioning,shapes,shadows,arrows,fit}
\usepackage{algpseudocode}

\usepackage{multicol}
\usepackage{graphicx}
\usepackage{array}
\usepackage{cite}

\newtheorem{Theoreme}{Theorem}[section]

\newtheorem{Proposition}[Theoreme]{Proposition}
\newtheorem{Lemme}[Theoreme]{Lemma}
\newtheorem{Definition}[Theoreme]{Definition}

\newtheorem{Remarque}[Theoreme]{Remark}

\DeclareMathOperator{\FQSym}{\mathbf{FQSym}}
\DeclareMathOperator{\PBT}{\mathbf{PBT}}
\newcommand{\FQSymm}{\FQSym^{(m)}} 

\newcommand{\PBTm}{\PBT^{(m)}} 

\newcommand{\Sym}[1]{\mathfrak{S}_{#1}} 
\newcommand{\Symmn}{\mathfrak{S}^{(m)}_{n}} 
\newcommand{\Symk}[2]{\mathfrak{S}^{(#2)}_{#1}} 

\DeclareMathOperator{\std}{std} 
\DeclareMathOperator{\Sup}{\vee} 
\DeclareMathOperator{\Inf}{\wedge} 
\DeclareMathOperator{\coinv}{coinv}
 
\DeclareMathOperator{\class}{\mathcal{C}} 
\DeclareMathOperator{\maxclass}{MaxClass} 
\DeclareMathOperator{\dt}{DT} 
\DeclareMathOperator{\Cinv}{\mathcal{I}} 

\newcommand{\Tamnm}{\mathcal{T}_{n}^{(m)}} 

\newcommand{\Metasmn}{\mathcal{MS}_n^{(m)}} 
\newcommand{\Metask}[2]{\mathcal{MS}_{#1}^{(#2)}} 
\newcommand{\Maxmn}{\mathcal{M}ax_n^{(m)}} 
\newcommand{\Maxk}[2]{\mathcal{M}ax_{#1}^{(#2)}} 
\newcommand{\Decmn}{\mathcal{DT}_n^{(m)}} 
\newcommand{\Decmmn}{\mathcal{DT}_n^{(m+1)}} 
\newcommand{\Deck}[2]{\mathcal{DT}_{#1}^{(#2)}} 
\newcommand{\MCmn}{\mathcal{MC}_n^{(m)}} 
\newcommand{\MCk}[2]{\mathcal{MC}_{#1}^{(#2)}} 
\newcommand{\sigch}[1]{\sigma^{(#1)}} 
\newcommand{\much}[1]{\mu^{(#1)}} 

\definecolor{darkGreen}{RGB}{23,103,1}
\newcommand{\red}[1]{\textbf{\textcolor{red}{#1}}}
\newcommand{\sred}[1]{\textcolor{red}{#1}}
\newcommand{\blue}[1]{\textcolor{blue}{#1}}
\newcommand{\green}[1]{\textcolor{darkGreen}{#1}}

\tikzstyle{Leaf} = [color = gray]

\tikzstyle{Red} = [color = red]
\tikzstyle{Blue} = [color = blue]
\tikzstyle{Green} = [color = darkGreen]
\tikzstyle{Gray} = [color = gray]

\tikzstyle{Path} = [line width = 1.2]
\tikzstyle{StrongPath} =  [line width=2.5]
\tikzstyle{DPoint} = [fill, radius=0.1]
\tikzstyle{Line1} = [dashed]
\tikzstyle{Line2} = [dotted, ultra thick]

\tikzstyle{Point} = [fill, radius=0.08]
\tikzstyle{RedPoint} = [color = red, fill, radius=0.08]
\tikzstyle{BluePoint} = [color = blue, fill, radius=0.08]
\tikzstyle{GreenPoint} = [color = darkGreen, fill, radius=0.08]
\tikzstyle{RedPath} = [color = red]
\tikzstyle{BluePath} = [color = blue]
\tikzstyle{GreenPath} = [color = darkGreen]
\tikzstyle{GrayPath} = [color = gray]
\tikzstyle{StrongPath} =  [line width=3.5]
\tikzstyle{StrongPath2} =  [line width=1.5]

\author{Viviane Pons}

\address{
Laboratoire de Recherche en Informatique \\
Bât 650 Ada Lovelace, Université Paris Sud, 91405 Orsay Cedex France
    }
    
\title[The metasylvester lattice]{A lattice on decreasing trees : the metasylvester lattice}


\begin{document}
\begin{abstract}
We introduce a new combinatorial structure: the metasylvester lattice on decreasing trees. It appears in the context of the $m$-Tamari lattices and other related $m$-generalizations. The metasylvester congruence has been recently introduced by Novelli and Thibon. We show that it defines a sublattice of the $m$-permutations where elements can be represented by decreasing labelled trees: the metasylvester lattice. We study the combinatorial properties of this new structure. In particular, we give different realizations of the lattice. The $m$-Tamari lattice is by definition a sublattice of our newly defined metasylvester lattice. It leads us to a new realization of the $m$-Tamari lattice, using certain chains of the classical Tamari lattice. 
\end{abstract}

\maketitle

\section{Introduction}
\label{sec:intro}

The Tamari lattice and its many generalizations have raised an increasing interest over the past few years. The lattice itself was introduced in \cite{Tamari1} by Tamari himself and is the keystone of many geometrical, algebraic and combinatorial constructions. The recent generalization to $m$-Tamari lattices by Bergeron and Préville-Ratelle \cite{BergmTamari} open many new combinatorial questions. The generalization of both the Malvenuto-Reutenauer and Loday-Ronco Hopf algebras \cite{MalReut, PBT1} to $m$-objects in the way of $\FQSym$ and $\PBT$ \cite{NCSF6, PBT2} was one of them which was studied in \cite{mTamAlg}.

It was shown in \cite{PBT2} that the Tamari lattice can be obtained from a congruence relation on permutations called the \emph{sylvester}\footnote{The word \emph{sylvester} is a litteral translation of the french word \emph{sylvestre} and, as so, is written in lower case.} congruence. in \cite{mTamAlg}, the same is done on a family of objects called $m$-permutations and then the $m$-Tamari lattice is obtained. The authors also introduce a new relation: the \emph{metasylvester} congruence. They show that the  metasylvester classes are in bijection with a certain family of labelled trees: the \emph{decreasing $(m+1)$-ary trees}. In the classical case where $m=1$, decreasing binary trees are actually in bijection with permutations themselves. And indeed, the metasylvester congruence is nontrivial only when $m>1$. In \cite{mTamAlg}, the authors are interested in the algebraic aspects: they describe the Hopf algebra derived from the congruence. It is shown to be both a subalgebra and quotient algebra of $\FQSymm$, the Hopf algebra of $m$-permutations. And $\PBTm$, the Hopf algebra on $(m+1)$-ary trees, is itself a subalgebra and quotient algebra of the metasylvester Hopf alegbra. 

In this paper, we study some combinatorial properties of the metasylvester classes. In particular, we introduce a new combinatorial structure: the metasylvester lattice on decreasing $(m+1)$-ary trees. It is a lattice \emph{in-between} the right weak order on $m$-permutations and the $m$-Tamari lattice. Indeed, it is a sublattice and a join-quotient lattice of the lattice on $m$-permutations and the $m$-Tamari lattice is a sublattice and a quotient lattice of it. Also, similarly to sylvester classes, metasylvester classes form intervals of the right weak order, and sylvester classes form intervals of the metasylvester lattice. Furthermore, we describe a very interesting realization of our new lattice in terms of chains of permutations. This leads us to a new realization of the $m$-Tamari lattice that we introduce at the end of the paper.

Section \ref{sub-sec:mperms} gives the definition of the lattice of $m$-permutation as an ideal of the right weak order. It also reminds some classical properties of the right weak order that are needed further on. The metaylvester lattice is defined in Section \ref{sub-sec:lattice-def} as a sublattice of the lattice on $m$-permutations. In Section \ref{sec:properties}, we explore the combinatorial properties of the lattice. We define the notion of \emph{tree-inversions} which is essential to our work and allows us to describe the cover relations. We give some important results in Section \ref{sub-sec:intervals}: metasylvester classes form intervals of the right weak order and the metasylvester lattice is a join-quotient of the $m$-permutations lattice. In Section \ref{sub-sec:chains}, we describe the realization of the lattice in terms of chains of permutations. Finally, in Section \ref{sub-sec:tam}, we summarize the links between the metasylvester and $m$-Tamari lattices and introduce a new realization of the latter.

\textbf{Notation:} in all examples, lattices and posets are represented with the \emph{smallest} element at the top of the lattice. 

\section{First definition of the metasylvester lattice}
\label{sec:metasylv}

\subsection{The lattice of $m$-permutations}
\label{sub-sec:mperms}

\begin{Definition}
An $m$-permutation of size $n$ is a permutation of the word $1^m 2^m \dots n^m$. We denote by $\Symmn$ the set of all $m$-permutations of size $n$.
\end{Definition}

\begin{figure}[ht]
\centering
\scalebox{0.8}{
\begin{tabular}{cc}

\def \wlev{1}
\def \hlev{1}

\begin{tikzpicture}
\node(T1) at (0,0){1122};
\node(T2) at (0,-1 * \hlev){1212};
\node(T3) at (-1 * \wlev, -2 * \hlev){2112};
\node(T4) at (\wlev, -2 * \hlev){1221};
\node(T5) at (0, -3 * \hlev){2121};
\node(T6) at (0, -4 * \hlev){2211};

\draw(T1) -- (T2);
\draw(T2) -- (T3);
\draw(T2) -- (T4);
\draw(T3) -- (T5);
\draw(T4) -- (T5);
\draw(T5) -- (T6);
\end{tikzpicture} &

\def \wlev{1}
\def \hlev{1}

\begin{tikzpicture}
\node(T1) at (0,0){1234};
\node(T2) at (0,-1 * \hlev){1324};
\node(T3) at (-1 * \wlev, -2 * \hlev){3124};
\node(T4) at (\wlev, -2 * \hlev){1342};
\node(T5) at (0, -3 * \hlev){3142};
\node(T6) at (0, -4 * \hlev){3412};

\draw(T1) -- (T2);
\draw(T2) -- (T3);
\draw(T2) -- (T4);
\draw(T3) -- (T5);
\draw(T4) -- (T5);
\draw(T5) -- (T6);
\end{tikzpicture}
\end{tabular}
}
\caption{Lattice on $\Symk{2}{2}$ as an ideal of the right weak lattice on $\Sym{4}$.}
\label{fig:mperms}
\end{figure}
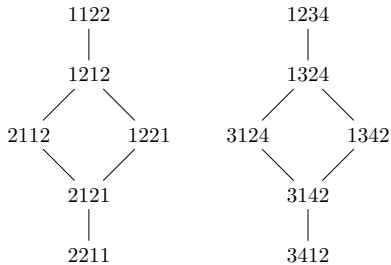

As an example, $122313$ is a $2$-permutation of size 3. These objects have been introduced before \cite{ME_THESE, mTamAlg}. They are used to generalize the Hopf algebra structures of $\FQSym$ and $\PBT$ to $m$-objects. In particular, the set of $m$-permutations of size $n$ naturally possess a lattice structure induced from the right weak order on permutations. Indeed, as explained in \cite{mTamAlg}, a $m$-permutation can be standardized into a classical permutation of size $m \times n$. The previous example gives $\std(122313) = 134526$. By this operation, $\Symmn$ forms an ideal of the right weak order on $\Sym{n}$ generated by $\std(n^m (n-1)^m \dots 1^m)$. An example is given Figure \ref{fig:mperms}: $\Symk{2}{2}$ is the ideal of $\Sym{4}$ generated by $3412$.

This structure on $m$-permutations is crucial for our work in this paper. We are going to use it to define the metasylvester lattice. First, let us recall some well-known properties of the right weak order and apply them to the $m$-permutations lattice. 

\begin{Definition}
A co-inversion of a permutation $\sigma \in \Sym{n}$ is a couple $(a,b)$ where $a < b$ and $b$ appears before $a$ in $\sigma$. We denote the set of co-inversions of $\sigma$ by $\coinv(\sigma)$.
\end{Definition}

As an example, $(2,4)$ is a co-inversion of $\sigma = 134526$. The permutation can be retrieved from the set of its co-inversions. It satisfies a few properties.

\begin{Proposition}
\label{prop:co-invs}
A set $L$ of couples $(a,b)$ with $1 \leq a \leq b \leq n$ corresponds to the co-inversions of a given permutation $\sigma$ if and only if
\begin{enumerate}
\item if $(a,b) \in L$ and $(b,c) \in L$, then $(a,c) \in L$,
\item if $(a,c) \in L$, then for all $b$ such that $a<b<c$, either $(a,b) \in L$ or $(b,c) \in L$.
\end{enumerate}
\end{Proposition}

The actual set of co-inversions is not needed to retrieve the permutation. Indeed, permutations are in bijection with their \emph{co-code}. The co-code of a permutation is a vector $v = v_1 \dots v_n$ where $v_i$ is the number of co-inversions of the form $(i,*)$ (it is the Lehmer code of the inverse of the permutation). As an example, the co-code of $23154$ is $20010$. A vector $v$ is a co-code if and only if $v_i \leq n-i$ for all $i$.

The notion of co-inversion still holds for $m$-permutations. Only, each letter is now repeated $m$ times. Thus, we call $a_i$ the $i^{th}$ $a$ of an $m$-permutation and a co-inversion $(a_i, b_j)$ means that the $j^{th}$ $b$ is before the $i^{th}$ $a$. As an example, the co-inversions of the $m$-permutations $2121$ are $(1_1, 2_1)$, $(1_2, 2_1)$, and $(1_2, 2_2)$. Note that if an $m$-permutation possesses the co-inversion $(a_i, b_j)$, it also possess all co-inversions $(a_k, b_j)$ for $k>i$. The following property is well-known for the right weak order and so still holds in the lattice of $m$-permutations.

\begin{Proposition}
For $\sigma, \mu \in \Symmn$, we have $\sigma \leq \mu$ if and only if   $\coinv(\sigma) \subseteq \coinv(\mu)$.
\end{Proposition}

The right weak order is a lattice, which means that for every two permutations $\sigma, \mu \in \Sym{n}$, the least upper bound (or join), $\sigma \Sup \mu$, and the greatest lower bound (or meet), $\sigma \Inf \mu$, of $\sigma$ and $\mu$ are well-defined. To compute $\sigma \Sup \mu$, one can read a minimal linear extension of a poset formed by the co-inversions of both $\sigma$ and $\mu$. The computation of $\sigma \Inf \mu$ is quite similar. As the order on $\Symmn$ forms an ideal of the right weak order, the meet and join operations are also well-defined on $m$-permutations.

\subsection{The metasylvester lattice}
\label{sub-sec:lattice-def}

In \cite{PBT2}, the authors define a monoid congruence on words called the \emph{sylvester} congruence. The permutations of $\Sym{n}$ can thus be divided into \emph{sylvester classes} which each forms an interval of the right weak order. It is then possible to define an order relation on these classes which is isomorphic to the Tamari lattice. As explained in \cite{ME_THESE} and \cite{mTamAlg}, the same can be done on $m$-permutations and we then obtain the $m$-Tamari lattice. In this paper, we are interested in a finer monoid congruence relation which was introduced in \cite{mTamAlg}, the \emph{metasylvester} congruence. It is defined by the transitive and reflexive closure of the relations

\begin{align}
\label{eq:metacongru-1}
ac\dots a &\equiv ca \dots a & (a < c),\\
\label{eq:metacongru-2}
b \dots ac \dots b &\equiv b \dots ca \dots b & (a < b < c).
\end{align}

As an example, $231312 \equiv 321132$ because we can do the following rewritings $231312 \rightarrow 321312 \rightarrow 321132$. For an $m$-permutation $\sigma$, we denote by $\class(\sigma)$ the class of $\sigma$, \emph{i.e.}, the set of all $m$-permutations that are equivalent to $\sigma$. In \cite{mTamAlg}, the authors use this congruence relation to define a sub-algebra of $\FQSymm$. In this paper, we are interested in the structure of the metasylvester classes in the $m$-permutations lattices.

We denote by $\Metasmn$ the set of metasylvester classes of $m$-permutations of size $n$. The set $\Metask{2}{2}$ consists of three classes: $\lbrace 1122 \rbrace$, $\lbrace 1212, 2112 \rbrace$, and  $\lbrace 1221, 2121, 2211 \rbrace$. More generally, here is the enumeration for different values of $n$ and $m$.

\begin{center}
\scalebox{0.8}{
\begin{tabular}{l|lllll}
$m\setminus n$  & 1 & 2 & 3  & 4    & 5     \\
\hline
1               & 1 & 2 & 6  & 24   & 120   \\
2               & 1 & 3 & 15 & 105  & 945   \\
3               & 1 & 4 & 28 & 280  & 3640  \\
4               & 1 & 5 & 45 & 585  & 9945  \\
5               & 1 & 6 & 66 & 1056 & 22176 \\

\end{tabular}
}
\end{center}

For $m=1$, the classes are reduced to one element and so correspond to the permutations. For a general $m$, the formula itself is quite simple,

\begin{equation}
\label{eq:nb-metasylv}
\mid \Metasmn \mid = (1+m)(1+2m)\dots (1+(n-1)m).
\end{equation}

It is proven in \cite{mTamAlg} by a bijection between $\Metasmn$ and $(m+1)$-ary decreasing trees. In particular, one can define a canonical element for each class.

\begin{Proposition}[Proposition 3.10 of \cite{mTamAlg}]
\label{prop:max}
For $\sigma \in \Symmn$, the class $\class(\sigma)$ contains a unique maximal element in the $m$-permutations lattice. It is the only element of the class which does not contain a sub-word $a\dots b \dots a$ with $a<b$. We denote it by $\maxclass(\class(\sigma))$, or directly $\maxclass(\sigma)$, and $\Maxmn$ is the set of all maximal elements of $\Metasmn$.
\end{Proposition}

As an example, $\Maxk{2}{2} = \lbrace 1122, 2112, 2211 \rbrace$.  The elements of $\Maxmn$ have an interesting recursive structure, they are \emph{Stirling permutations} \cite{Stirling} and are in bijection with \emph{decreasing $(m+1)$-ary trees}.

\begin{Definition}
We call decreasing $m$-ary tree, a planar $m$-ary tree whose internal nodes are labelled with distinct labels $1, \dots, n$ such that a node label is always greater than all labels of its descendent nodes. 

We denote the set of all $m$-ary decreasing trees by $\Decmn$.
\end{Definition}

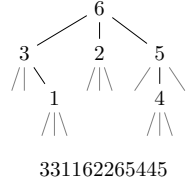
\begin{figure}[ht]
\begin{center}
\scalebox{0.8}{
\begin{tabular}{c}
\begin{tikzpicture}
\node (T0) at (11.250, 0.000){6};
\node (T00) at (10.000, -0.750){3};
\node(T000) at (9.750, -1.500){};
\node(T001) at (10.000, -1.500){};
\node (T002) at (10.500, -1.500){1};
\node(T0020) at (10.250, -2.250){};
\node(T0021) at (10.500, -2.250){};
\node(T0022) at (10.750, -2.250){};
\draw[Leaf] (T002) -- (T0020);
\draw[Leaf] (T002) -- (T0021);
\draw[Leaf] (T002) -- (T0022);
\draw[Leaf] (T00) -- (T000);
\draw[Leaf] (T00) -- (T001);
\draw (T00) -- (T002);
\node (T01) at (11.250, -0.750){2};
\node(T010) at (11.000, -1.500){};
\node(T011) at (11.250, -1.500){};
\node(T012) at (11.500, -1.500){};
\draw[Leaf] (T01) -- (T010);
\draw[Leaf] (T01) -- (T011);
\draw[Leaf] (T01) -- (T012);
\node (T02) at (12.250, -0.750){5};
\node(T020) at (11.750, -1.500){};
\node (T021) at (12.250, -1.500){4};
\node(T0210) at (12.000, -2.250){};
\node(T0211) at (12.250, -2.250){};
\node(T0212) at (12.500, -2.250){};
\draw[Leaf] (T021) -- (T0210);
\draw[Leaf] (T021) -- (T0211);
\draw[Leaf] (T021) -- (T0212);
\node(T022) at (12.750, -1.500){};
\draw[Leaf] (T02) -- (T020);
\draw (T02) -- (T021);
\draw[Leaf] (T02) -- (T022);
\draw (T0) -- (T00);
\draw (T0) -- (T01);
\draw (T0) -- (T02);
\end{tikzpicture}
\\
331162265445
\end{tabular}
}
\end{center}
\caption{An element of $\Deck{6}{3}$ and its corresponding element of $\Maxk{6}{2}$.}
\label{fig:decreasing-tree}
\end{figure}

To obtain a decreasing tree from an element $\sigma \in \Maxmn$, one follows a recursive algorithm which generalizes the decreasing binary tree of a permutation \cite{PBT2}. The root of the tree is always labelled by $n$, the maximal letter of $\sigma$. Then we use the occurrences of $n$ to divide $\sigma$ into $m+1$ factors. As an example the $2$-permutation 3311\textbf{6}22\textbf{6}5445 is divided into $3311$, $22$, $5445$. As $\sigma$ avoids $a \dots b \dots a$, all letters of same value appear in the same block. Those blocks correspond to the $m+1$ subtrees and we apply the procedure recursively. The inverse operation can also be easily described. The $m$-permutation corresponds to a recursive traversal of the decreasing tree: first subtree on the left, root, second subtree, root, \dots, root, last subtree on the right. We denote by $\maxclass(T)$ the $m$-permutation of a decreasing tree $T$. An example is given in Figure \ref{fig:decreasing-tree}.

We denote by $\dt : \Maxmn \rightarrow \Decmmn$ the function giving the decreasing tree of a maximal class element. This function can actually be extended to all $m$-permutations by $\dt(\sigma) = \dt(\maxclass(\sigma))$ and then becomes a surjection from $\Symmn$ to $\Decmmn$. A direct algorithm to construct $\dt(\sigma)$ is given in \cite[Algorithm 3.8]{mTamAlg}. We propose a new description of it in the next Section.

The maximal class elements form a subset of the $m$-permutations, thus we can define a partial order induced by the order on $m$-permutations. We call this partial order the \emph{metasylvester poset}.

\begin{Theoreme}
The metasylvester poset is a lattice.
\end{Theoreme}

This is proved by showing that $\Maxmn$ is stable through join and meet operations. An example of the lattice for $n=3$ and $m=2$ is given in Figure \ref{fig:metasylv-lattice}.

\begin{figure}[ht]
\scalebox{0.7}{
\begin{tabular}{ccc}
\begin{minipage}{3in}

\def \wlev{1.5}
\def \hlev{1.5}

\begin{tikzpicture}
\node(P112233) at (0,0){112233};
\node(P211233) at (-1 * \wlev, -1 * \hlev){211233};
\node(P113223) at ( 1 * \wlev, -1 * \hlev){113223};

\node(P221133) at (-2 * \wlev, -2 * \hlev){221133};
\node(P311223) at ( 1 * \wlev, -2 * \hlev){311223};
\node(P113322) at ( 2 * \wlev, -2 * \hlev){113322};

\node(P223113) at (-2 * \wlev, -3 * \hlev){223113};
\node(P321123) at ( 0 * \wlev, -3 * \hlev){321123};
\node(P311322) at ( 2 * \wlev, -3 * \hlev){311322};

\node(P223311) at (-2 * \wlev, -4 * \hlev){223311};
\node(P322113) at (-1 * \wlev, -4 * \hlev){322113};
\node(P331122) at ( 2 * \wlev, -4 * \hlev){331122};

\node(P322311) at (-1 * \wlev, -5 * \hlev){322311};
\node(P332112) at ( 1 * \wlev, -5 * \hlev){332112};

\node(P332211) at ( 0 * \wlev, -6 * \hlev){332211};

\draw (P112233) -- (P211233);
\draw (P112233) -- (P113223);

\draw(P211233) -- (P221133);
\draw(P211233) -- (P321123);
\draw(P113223) -- (P311223);
\draw(P113223) -- (P113322);

\draw(P221133) -- (P223113);
\draw(P311223) -- (P321123);
\draw(P311223) -- (P311322);
\draw(P113322) -- (P311322);

\draw(P223113) -- (P223311);
\draw(P223113) -- (P322113);
\draw(P321123) -- (P322113);
\draw(P321123) -- (P332112);
\draw(P311322) -- (P331122);

\draw(P223311) -- (P322311);
\draw(P322113) -- (P322311);
\draw(P331122) -- (P332112);

\draw(P322311) -- (P332211);
\draw(P332112) -- (P332211);
\end{tikzpicture}
\end{minipage}
&
\begin{minipage}{3in}

\def \wlev{1.6}
\def \hlev{1.6}
\def \fpath{figures/decreasing_trees/}
\def \sscale{0.5}

\begin{tikzpicture}
\node(P112233) at (0,0){\scalebox{\sscale}{\input{\fpath DT3-3-a}}};
\node(P211233) at (-1 * \wlev, -1 * \hlev){\scalebox{\sscale}{\input{\fpath DT3-3-b}}};
\node(P113223) at ( 1 * \wlev, -1 * \hlev){\scalebox{\sscale}{\input{\fpath DT3-3-f}}};

\node(P221133) at (-2 * \wlev, -2 * \hlev){\scalebox{\sscale}{\input{\fpath DT3-3-c}}};
\node(P311223) at ( 1 * \wlev, -2 * \hlev){\scalebox{\sscale}{\input{\fpath DT3-3-g}}};
\node(P113322) at ( 2 * \wlev, -2 * \hlev){\scalebox{\sscale}{\input{\fpath DT3-3-k}}};

\node(P223113) at (-2 * \wlev, -3 * \hlev){\scalebox{\sscale}{\input{\fpath DT3-3-d}}};
\node(P321123) at ( 0 * \wlev, -3 * \hlev){\scalebox{\sscale}{\input{\fpath DT3-3-h}}};
\node(P311322) at ( 2 * \wlev, -3 * \hlev){\scalebox{\sscale}{\input{\fpath DT3-3-l}}};

\node(P223311) at (-2 * \wlev, -4 * \hlev){\scalebox{\sscale}{\input{\fpath DT3-3-e}}};
\node(P322113) at (-1 * \wlev, -4 * \hlev){\scalebox{\sscale}{\input{\fpath DT3-3-i}}};
\node(P331122) at ( 2 * \wlev, -4 * \hlev){\scalebox{\sscale}{\input{\fpath DT3-3-m}}};

\node(P322311) at (-1 * \wlev, -5 * \hlev){\scalebox{\sscale}{\input{\fpath DT3-3-j}}};
\node(P332112) at ( 1 * \wlev, -5 * \hlev){\scalebox{\sscale}{\input{\fpath DT3-3-n}}};

\node(P332211) at ( 0 * \wlev, -6 * \hlev){\scalebox{\sscale}{\input{\fpath DT3-3-o}}};

\draw (P112233) -- (P211233);
\draw (P112233) -- (P113223);

\draw(P211233) -- (P221133);
\draw(P211233) -- (P321123);
\draw(P113223) -- (P311223);
\draw(P113223) -- (P113322);

\draw(P221133) -- (P223113);
\draw(P311223) -- (P321123);
\draw(P311223) -- (P311322);
\draw(P113322) -- (P311322);

\draw(P223113) -- (P223311);
\draw(P223113) -- (P322113);
\draw(P321123) -- (P322113);
\draw(P321123) -- (P332112);
\draw(P311322) -- (P331122);

\draw(P223311) -- (P322311);
\draw(P322113) -- (P322311);
\draw(P331122) -- (P332112);

\draw(P322311) -- (P332211);
\draw(P332112) -- (P332211);
\end{tikzpicture}
\end{minipage}
&
\begin{minipage}{3in}
\def \wlev{1.5}
\def \hlev{1.5}

\begin{tikzpicture}[every text node part/.style={align=center}]

\node(P123_123) at (0,0){
123 \\
123
};

\node(P123_213) at (-1 * \wlev, -1 * \hlev){
123 \\
213
};

\node(P123_132) at (1 * \wlev, -1 * \hlev){
123 \\
132
};

\node(P213_213) at (-2 * \wlev, -2 * \hlev){
213 \\
213
};

\node(P123_312) at (1 * \wlev, -2 * \hlev){
123 \\
312
};

\node(P132_132) at (2 * \wlev, -2 * \hlev){
132 \\
132
};

\node(P213_231) at (-2 * \wlev, -3 * \hlev){
213 \\
231
};

\node(P123_321) at (0 * \wlev, -3 * \hlev){
123 \\
321
};

\node(P132_312) at (2 * \wlev, -3 * \hlev){
132 \\
312
};

\node(P231_231) at (-2 * \wlev, -4 * \hlev){
231 \\
231
};

\node(P213_321) at (-1 * \wlev, -4 * \hlev){
213 \\
321
};

\node(P312_312) at (2 * \wlev, -4 * \hlev){
312 \\
312
};

\node(P231_321) at (-1 * \wlev, -5 * \hlev){
231 \\
321
};

\node(P312_321) at (1 * \wlev, -5 * \hlev){
312 \\
321
};

\node(P321_321) at (0 * \wlev, -6 * \hlev){
321 \\
321
};

\draw (P123_123) -- (P123_213);
\draw (P123_123) -- (P123_132);
\draw (P123_213) -- (P213_213);
\draw (P123_213) -- (P123_321);
\draw (P123_132) -- (P123_312);
\draw (P123_132) -- (P132_132);
\draw (P213_213) -- (P213_231);
\draw (P123_312) -- (P123_321);
\draw (P123_312) -- (P132_312);
\draw (P132_132) -- (P132_312);
\draw (P213_231) -- (P231_231);
\draw (P213_231) -- (P213_321);
\draw (P123_321) -- (P213_321);
\draw (P123_321) -- (P312_321);
\draw (P132_312) -- (P312_312);
\draw (P231_231) -- (P231_321);
\draw (P213_321) -- (P231_321);
\draw (P312_312) -- (P312_321);
\draw (P231_321) -- (P321_321);
\draw (P312_321) -- (P321_321);

\end{tikzpicture}
\end{minipage}
\end{tabular}
}
\caption{The metasylvester lattice on $\Maxk{3}{2}$, $\Deck{3}{3}$, and $\MCk{3}{2}$.}
\label{fig:metasylv-lattice}
\end{figure}

\section{Cover relations and combinatorial properties}
\label{sec:properties}

\subsection{Tree-inversions}
\label{sub-sec:tree-invs}

One purpose of this paper is to give different realizations of our newly defined metasylvester lattice. These realizations are given in terms of different combinatorial objects (we have seen already maximal class elements $\Maxmn$ and decreasing trees $\Decmn$) and the key to all bijections lies into the notion of \emph{tree-inversion} that we define now. It will also allow us to describe the cover relations of the lattice and many other properties of metasylvester classes. 

\begin{Definition}
\label{def:tree-inversion}
Let $\sigma$ be an $m$-permutation of size $n$ and $a,b$ be such that $1 \leq a < b \leq n$. We say that  $\sigma$ contains the \emph{tree-inversion} $(a,b_i)$ with $1 \leq i \leq m$, if there is a co-inversion $(a,b_i)$ in $\maxclass(\sigma)$, \emph{i.e.}, all letters $a$ are placed after the $i^{th}$ $b$ (remind that either all letters $a$ or none of them can be placed after a letter $b$).

If $(a,b_i)$ is not a tree-inversion, we say that it is \emph{tree-sorted}.
\end{Definition}

If $\sigma$ is already a maximal class element, then its number of tree-inversions is exactly its number of co-inversions divided by $m$ (the co-inversions $(a_1, b_i), \dots, (a_m, b_i)$ only contribute to 1 tree-inversion). It can also be read directly on the decreasing tree. A tree admits a tree-inversion $(a,b_i)$ if $a<b$ and: either $a$ belongs to the $j^{th}$ subtree of $b$ with $j > i$, or $a$ and $b$ have a common ancestor and $a$ is in a righter subtree than $b$. Equivalently, $(a,b_i)$ is tree-sorted if either $a$ belongs to the $j^{th}$ subtree of $b$ with $j \leq i$ or $a$ and $b$ have a common ancestor and $a$ is in a lefter subtree than $b$. As an example, in Figure \ref{fig:decreasing-tree}, $(4,5_1)$ is a tree-inversion because $4$ is in the second subtree of $5$, and $(2,3_1)$ and $(2,3_2)$ are tree-inversions because $2$ is on the right of $3$. On the other hand $(1,2_1)$ and $(1,2_2)$ are tree-sorted because $1$ is on the left of $2$. 

From the definition, the co-inversions of a maximal class element can be directly retrieved from the tree-inversions: a tree-inversion $(a,b_i)$ gives all co-inversions $(a_j,b_i)$. Thus, a set of tree-inversions is a unique identifier of a metasylvester class. More precisely, we can describe the lists we obtain.

\begin{Proposition}
\label{prop:valid-inv-list}
A set $S$ of couples $(a,b_i)$ with $1 \leq a \leq b \leq n$ and $1 \leq i \leq m$ is a \emph{valid} set of tree-inversions, \emph{i.e.}, there is a $m$-permutation $\sigma \in \Maxmn$ such that the tree-inversions of $\sigma$ are exactly $S$, if and only if
\begin{enumerate}
\item if $(a,b_j) \in S$, then $(a,b_i) \in S$ for all $i \leq j$,
\label{cond:valid-inv-list-m}
\item if $(a,b_j) \in S$ and $(b,c_i) \in S$, then $(a,c_i) \in S$,
\label{cond:valid-inv-list-transitivity}
\item if $(a,c_i) \in S$, then, for all $b$ such that $a<b<c$, either $(a,b_j) \in S$ for all $j$, or $(b,c_i) \in S$.
\label{cond:valid-inv-list-inv}
\end{enumerate}
\end{Proposition}

This proposition is just a direct consequence of the properties of co-inversions of the right weak order given in Proposition \ref{prop:co-invs}. Indeed, Conditions \eqref{cond:valid-inv-list-transitivity} and \eqref{cond:valid-inv-list-inv} are the direct transcription of Proposition~\ref{prop:co-invs} in the case of elements of $\Maxmn$. And Condition \eqref{cond:valid-inv-list-m} is from the specificity of $m$-permutations. Moreover, we can define the \emph{tree-code} of a metasylvester class the same way we defined the co-code. If $\sigma \in \Maxmn$, the co-inversions $(a_i,*)$ are the same for all $1 \leq i \leq m$, and then, the co-code is of the form $v_1 \dots v_1 v_2 \dots v_2 \dots v_n \dots v_n$. We define the tree-code by $v_1 v_2 \dots v_n$ where $v_i$ is the number of tree-inversions $(i,*)$. As an example, the maximal class element $\sigma = 331162265445$ has a co-code of $223300332200$ and the tree-code of its class is then $230320$. One can check that there are indeed two tree-inversions with 1, $(1,3_1)$ and $(1,3_2)$, three tree-inversions with 2, $(2,3_1), (2,3_2), (2,6_1)$ and so on. A vector $v=v_1 \dots v_n$ is a tree-code if and only if $0 \leq v_i \leq (n-i) m +1$ for all $1 \leq i \leq n$. It is clear that it is necessary, and it is indeed sufficient because we can get back a valid co-code from it. In particular, this gives a direct bijective proof of \eqref{eq:nb-metasylv}.

We can read tree-inversions from both maximal class elements and decreasing trees. The following proposition allows to identify directly the tree-inversions of any given $m$-permutation. 

\begin{Proposition}
\label{prop:tree-sorted}
Let $\sigma$ be an $m$-permutation, then $(a,b_i)$ is tree-sorted if and only if 
\begin{enumerate}
\item the last occurrence of $a$ is placed before the $i^{th}$ occurrence of $b$,
\label{cond:tree-sorted-1}
\item for all $a'$ such that $a < a' < b$ and there is an occurrence of $a'$ placed before the last occurrence $a$, $(a',b_i)$ is also tree-sorted.
\label{cond:tree-sorted-2}
\end{enumerate}

\end{Proposition}

As an example, in $\sigma = 12132434$, $(2,4_1)$ is not tree-sorted because $(3,4_1)$ is not. By consequence, $(1,4_1)$ is also not tree-sorted. The complete set of tree-inversions gives: $\lbrace(1,2_1), (1,3_1), (2,3_1), (1,4_1)$, $(2,4_1)$, $(3,4_1)\rbrace$. The proof goes in two steps: first notice that this new definition is indeed consistent with Definition \ref{def:tree-inversion} for maximal class elements, then show that the list keeps unchanged through elementary rewritings \eqref{eq:metacongru-1} and \eqref{eq:metacongru-2}. From this proposition, we get that the tree-inversions of an $m$-permutation $\sigma$ are the transitive closure through Condition \eqref{cond:valid-inv-list-transitivity} of Proposition \ref{prop:valid-inv-list} of the set of co-inversions $(a_m,b_i)$ of $\sigma$ with $1 \leq a \leq b \leq n$ and $1 \leq i \leq m$. Indeed, Proposition \ref{prop:valid-inv-list} can be read:  $(a,b_i)$ is a tree-inversion of $\sigma$ if and only if $(a_m, b_i)$ is a co-inversion of $\sigma$ or if there is $a' > a$ with $(a_m,a'_j) \in \coinv(\sigma)$ and $(a',b_i)$ is a tree-inversion of $\sigma$. Note that only one condition of Proposition \ref{prop:valid-inv-list} is needed to close the set: Conditions \eqref{cond:valid-inv-list-m} and \eqref{cond:valid-inv-list-inv} are already satisfied by the set of co-inversions $(a_m,b_i)$ of $\sigma$.

As we are now able to get the set of tree-inversions of any $m$-permutation $\sigma$, we can easily construct its decreasing tree. Indeed, the root of the tree is always labelled $n$, the maximal letter of the permutation. The tree-inversions involving $n$ give us the composition of the subtrees which we can construct recursively. Note that a similar algorithm was given in \cite[Algorithm 3.8]{mTamAlg} without the notion of tree-inversions. As an example if $\sigma = 133126245465$, the tree inversions $(*,6_2)$ are $(5,6_2)$ and $(4,6_2)$ which means that $5$ and $4$ form the rightest subtree. In the same way, $2$ is in the second subtree because there is a tree-inversion $(2,6_1)$. We then construct the decreasing trees of $1331$, $22$, and $4545$, and we obtain the tree of Figure \ref{fig:decreasing-tree}. It can be checked that $\sigma$ is indeed equivalent to the maximal element obtained from the tree.

\begin{Lemme}
\label{lem:tree-inversions-inclusion}
Let $\sigma, \mu \in \Symmn$ such that $\sigma \leq \mu$ then the tree-inversions of $\sigma$ are included in the tree-inversions of $\mu$. The inclusion is a strict one if and only if $\sigma \not\equiv \mu$
\end{Lemme}

Note that there is no reciprocity in general: the tree-inversions of $131223$ are included in the tree-inversions of $121332$ but the two permutations are not comparable in the $m$-permutation lattice. Nevertheless, there is reciprocity for maximal class elements and the \emph{classes} of the two permutations are indeed comparable in the metasylvester lattice.

\begin{Lemme}
\label{lem:interval-closure}
If $\sigma \leq \mu$ and $\sigma \equiv \mu$ then $\sigma \equiv \nu$ for all $\nu$ such that $\sigma \leq \nu \leq \mu$.
\end{Lemme}

Lemma \ref{lem:interval-closure} is a direct consequence of Lemma \ref{lem:tree-inversions-inclusion}. We use it as well as some specific properties of the $m$-permutations to describe the cover relations of the metasylvester lattice.

\begin{Proposition}
Let $\sigma \in \Maxmn$, then the successors of $\sigma$ in the metasylvester lattice are all permutations $\mu = \maxclass(\sigma')$ such that $\sigma'$ is a successor of $\sigma$ in the $m$-permutations lattice. They are all $m$-permutations $\mu$ such that $a < b$ and $\sigma = u~a_1~v~a_m~b~w$ and $\mu = u~b~a_1~v~a_m~w$ where $u,v$, and $w$ are words. (Note that because $\sigma \in \Maxmn$, then all letters of $v$ are smaller or equal to $a$).
\end{Proposition}

As an example, the successors of $22311344$ are $32211344$, $22331144$, and $22431134$.

\subsection{Interval and semi-quotient properties}
\label{sub-sec:intervals}

We said already that the metasylvester relation is a finer version of the more classical sylvester relation which is related to the Tamari and $m$-Tamari lattices. In particular, a sylvester class is formed by the union of many metasylvester classes. In this section, we want to investigate which properties of the sylvester classes are still satisfied by the metasylvester classes. A well-known property of the syvester classes is that they form intervals of the right weak order \cite{PBT2}. It is also true for metasylvester classes.

\begin{Proposition}
\label{prop:interval}
Let $C$ be metasylvester class, \emph{i.e.}, a set of permutations that are equivalent to one another, then $C$ forms an interval of the $m$-permutations lattice.
\end{Proposition}

From Lemma \ref{lem:interval-closure}, we already knew that the class was closed by interval. Also, it has a unique maximal element. We prove that there is a unique minimal element using \emph{Newman lemma} (or diamond lemma) on rewriting rules. Note that we prove the uniqueness of the minimal element but do not give the element itself. It is possible to describe the algorithm giving the minimal element from the decreasing tree but we could not find any interesting properties of the set of minimal elements. Indeed, unlike the sylvester case, the sub-poset formed by the minimal class elements is not the same as the metasylvester lattice and is not a lattice.

In the case of the classical sylvester relation, we know that the Tamari lattice can be seen either as a sub-lattice of the right weak order (by taking the subset of maximal class elements) or as quotient lattice when quotienting by the congruence relation. We just proved that the first property is also true for the metasylvester lattice but the second one is not. Indeed, for $\sigma = 121332$ and $\mu=131223$, we have

\begin{align}
\maxclass(\sigma \Inf \mu) &= \maxclass(112323) = 113223, \\
\maxclass(\sigma) \Inf \maxclass(\mu) &= 332112 \Inf 311223 = 311223.
\end{align}

And so $\class(\sigma \Inf \mu) \neq \class(\sigma) \Inf \class(\mu)$. Nevertheless, we still have a weaker quotient property and we can understand the metasylvester lattice as being a semi-quotient lattice of the $m$-permutation lattice for the join operation. This is a direct consequence of the properties of tree-inversions.

\begin{Proposition}
\label{prop:semi-quotient}
Let $\sigma$ and $\mu$ in $\Symmn$, then
\begin{align}
\class(\sigma \Sup \mu) = \class(\sigma) \Sup \class(\mu).
\end{align}

\end{Proposition}

\section{Other realization and $m$-Tamari}
\label{sec:other}

\subsection{Chains of permutations}
\label{sub-sec:chains}

The set of tree-inversions is the key object when working on the metasylvester lattice. We have seen that  it can be interpreted either as an $m$-permutation or as a decreasing tree. We now introduce a third combinatorial object, in bijection with the previous ones, that also derives from the tree-inversions.

\begin{Definition}
\label{def:metachain}
A \emph{metasylvester $m$-chain} $c = (\sigch{m}, \sigch{m-1}, \dots, \sigch{1})$ of size $n$ is a list of $m$-permutations which satisfy
\begin{equation}
\sigma^{(m)} \leq \sigma^{(m-1)} \leq \dots \leq \sigma^{(1)} 
\end{equation}
for the right weak order and such that, for all $i < j$, the permutation $(\sigma^{(j)})^{-1}\sigma^{(i)}$ avoids the pattern $231$.

We denote the set of metasylvester $m$-chains of size $n$ by $\MCmn$. 
\end{Definition}

As an example, the metasylvester 2-chains of size 3 are exactly all couples $\sigma \leq \mu$ except for $(123,231)$ and $(132,321)$, which gives 15 chains. By a simple introspection, one conjectures that the size of $\MCmn$ is given by \eqref{eq:nb-metasylv}. The purpose of this section is to give an explicit bijection between metasylvester $m$-chains and metasylvester classes.

\begin{Remarque}
\label{rem:avoid}
For two permutations $\sigma \leq \mu$, we have that $\sigma^{-1} \mu$ contains the pattern $231$ if and only if either $\sigma$ contains the subword $a \dots b \dots c$ and $\mu$ the subword $b \dots c \dots a$ or $\sigma$ contains the subword $a \dots c \dots b$ and $\mu$ the subword $c \dots b \dots a$ with $a < b < c$. 
\end{Remarque}

\begin{Proposition}
\label{prop:chains-tree-invs}
Let $c = (\sigch{m}, \dots, \sigch{1})$ be a metasylvester $m$-chain. We define the set $\Cinv(c)$ by

\begin{equation}
\Cinv(c) = \lbrace (a,b_i); (a,b) \in \coinv(\sigch{i}) \rbrace.
\end{equation} 

Then $\Cinv(c)$ is the set of tree-inversions of a given metasylvester class.
\end{Proposition}

As an example, if $c = (213, 231)$ then $\Cinv(c) = \lbrace (1,2_1), (1,3_1), (1, 2_2) \rbrace$ is the set of tree-inversions of the class of $223113$. We use Proposition \ref{prop:valid-inv-list} to prove that $\Cinv(c)$ is indeed a set of tree-inversions. Condition~\eqref{cond:valid-inv-list-m} comes from the fact that $\sigch{m} \leq \dots \leq \sigch{1}$. To prove Conditions \eqref{cond:valid-inv-list-transitivity} and \eqref{cond:valid-inv-list-inv}, we need to use the pattern avoidance property of metasylvester $m$-chains. The same kind of argument allows us to prove a reciprocal property.

\begin{Proposition}
\label{prop:tree-invs-chains}
Let $\mathcal{C}$ be a metasylvester class of $\Metasmn$ and $\mathcal{I}$ its set of tree-inversions. For all $1 \leq i \leq m$, we define $s_i = \lbrace (a,b); (a,b_i) \in \mathcal{I} \rbrace$. Then for each $s_i$, there is a permutation $\sigch{i}$ such that $\coinv(\sigch{i}) = s_i$ and $c = (\sigch{m}, \sigch{m-1}, \dots, \sigch{1})$ is a metasylvester $m$-chain. 
\end{Proposition}

\begin{Theoreme}
\label{thm:metasylv-chains}
Let $\psi : \MCmn \rightarrow \Metasmn$ the function that associates a metasylvester class to a metasylvester $m$-chain as described in Proposition \ref{prop:chains-tree-invs}. Then $\psi$ is a bijection between $\MCmn$ and $\Metasmn$.
\end{Theoreme}

Indeed, Proposition \ref{prop:chains-tree-invs} assures us that $\psi$ is well-defined. It it quite clear that it is injective and Proposition \ref{prop:tree-invs-chains} give us its inverse which makes it bijective. The bijection allows us to realize the metasylvester lattice on metasylvester chains (see the third realization in Figure \ref{fig:metasylv-lattice}). The order relation is immediate from the bijection itself. Indeed, the chains can be obtained from the tree-inversions list which actually corresponds to the co-inversions of the permutations. 

\begin{Proposition}
Let $c_1 = (\sigch{m}, \dots, \sigch{1})$ and $c_2 = (\much{m}, \dots, \much{1})$ , then $c_1 \leq c_2$ in the metasylvester lattice if and only if $\sigch{i} \leq \much{i}$ for all $1 \leq i \leq m$ in the right weak order.
\end{Proposition}

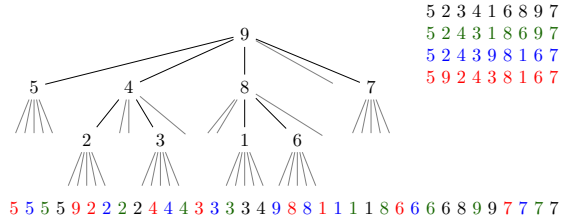
\begin{figure}[ht]
\begin{center}
\scalebox{0.7}{
\newcommand{\oa}[1]{\sred{#1} }
\newcommand{\ob}[1]{\blue{#1} }
\newcommand{\oc}[1]{\green{#1} }
\newcommand{\od}[1]{#1 }

\newcommand{\allof}[1]{\oa{#1}\ob{#1}\oc{#1}\od{#1}}

\begin{tabular}{cc}
\begin{tikzpicture}[baseline=-0.3cm]
\node (T0) at (4.400, 0.000){9};
\node (T00) at (0.400, -1.000){5};
\node(T000) at (0.000, -2.000){};
\node(T001) at (0.200, -2.000){};
\node(T002) at (0.400, -2.000){};
\node(T003) at (0.600, -2.000){};
\node(T004) at (0.800, -2.000){};
\draw[Leaf] (T00) -- (T000);
\draw[Leaf] (T00) -- (T001);
\draw[Leaf] (T00) -- (T002);
\draw[Leaf] (T00) -- (T003);
\draw[Leaf] (T00) -- (T004);
\node (T01) at (2.200, -1.000){4};
\node (T010) at (1.400, -2.000){2};
\node(T0100) at (1.000, -3.000){};
\node(T0101) at (1.200, -3.000){};
\node(T0102) at (1.400, -3.000){};
\node(T0103) at (1.600, -3.000){};
\node(T0104) at (1.800, -3.000){};
\draw[Leaf] (T010) -- (T0100);
\draw[Leaf] (T010) -- (T0101);
\draw[Leaf] (T010) -- (T0102);
\draw[Leaf] (T010) -- (T0103);
\draw[Leaf] (T010) -- (T0104);
\node(T011) at (2.000, -2.000){};
\node(T012) at (2.200, -2.000){};
\node (T013) at (2.800, -2.000){3};
\node(T0130) at (2.400, -3.000){};
\node(T0131) at (2.600, -3.000){};
\node(T0132) at (2.800, -3.000){};
\node(T0133) at (3.000, -3.000){};
\node(T0134) at (3.200, -3.000){};
\draw[Leaf] (T013) -- (T0130);
\draw[Leaf] (T013) -- (T0131);
\draw[Leaf] (T013) -- (T0132);
\draw[Leaf] (T013) -- (T0133);
\draw[Leaf] (T013) -- (T0134);
\node(T014) at (3.400, -2.000){};
\draw (T01) -- (T010);
\draw[Leaf] (T01) -- (T011);
\draw[Leaf] (T01) -- (T012);
\draw (T01) -- (T013);
\draw[Leaf] (T01) -- (T014);
\node (T02) at (4.400, -1.000){8};
\node(T020) at (3.600, -2.000){};
\node(T021) at (3.800, -2.000){};
\node (T022) at (4.400, -2.000){1};
\node(T0220) at (4.000, -3.000){};
\node(T0221) at (4.200, -3.000){};
\node(T0222) at (4.400, -3.000){};
\node(T0223) at (4.600, -3.000){};
\node(T0224) at (4.800, -3.000){};
\draw[Leaf] (T022) -- (T0220);
\draw[Leaf] (T022) -- (T0221);
\draw[Leaf] (T022) -- (T0222);
\draw[Leaf] (T022) -- (T0223);
\draw[Leaf] (T022) -- (T0224);
\node (T023) at (5.400, -2.000){6};
\node(T0230) at (5.000, -3.000){};
\node(T0231) at (5.200, -3.000){};
\node(T0232) at (5.400, -3.000){};
\node(T0233) at (5.600, -3.000){};
\node(T0234) at (5.800, -3.000){};
\draw[Leaf] (T023) -- (T0230);
\draw[Leaf] (T023) -- (T0231);
\draw[Leaf] (T023) -- (T0232);
\draw[Leaf] (T023) -- (T0233);
\draw[Leaf] (T023) -- (T0234);
\node(T024) at (6.000, -2.000){};
\draw[Leaf] (T02) -- (T020);
\draw[Leaf] (T02) -- (T021);
\draw (T02) -- (T022);
\draw (T02) -- (T023);
\draw[Leaf] (T02) -- (T024);
\node(T03) at (6.200, -1.000){};
\node (T04) at (6.800, -1.000){7};
\node(T040) at (6.400, -2.000){};
\node(T041) at (6.600, -2.000){};
\node(T042) at (6.800, -2.000){};
\node(T043) at (7.000, -2.000){};
\node(T044) at (7.200, -2.000){};
\draw[Leaf] (T04) -- (T040);
\draw[Leaf] (T04) -- (T041);
\draw[Leaf] (T04) -- (T042);
\draw[Leaf] (T04) -- (T043);
\draw[Leaf] (T04) -- (T044);
\draw (T0) -- (T00);
\draw (T0) -- (T01);
\draw (T0) -- (T02);
\draw[Leaf] (T0) -- (T03);
\draw (T0) -- (T04);
\end{tikzpicture}
&
\begin{tabular}{c}
5 2 3 4 1 6 8 9 7 \\
\green{5 2 4 3 1 8 6 9 7} \\
\blue{5 2 4 3 9 8 1 6 7} \\
\sred{5 9 2 4 3 8 1 6 7}
\end{tabular}
\\
\multicolumn{2}{c}{
\allof{5}\oa{9}\allof{2}\oa{4}\ob{4}\oc{4}\allof{3}\od{4}\ob{9}\oa{8}\ob{8}\allof{1}\oc{8}\allof{6}\od{8}\oc{9}\od{9}\allof{7}

}
\end{tabular}
}
\end{center}
\caption{A 5-decreasing tree, its corresponding metasylvester 4-chain and maximal class element} 
\label{fig:decreasing-tree-bijection}
\end{figure}

The bijection between tree-inversions lists and metasylvester chains can be nicely described directly on decreasing trees. Let $T$ be a $(m+1)$-ary decreasing tree whose subtrees are $T_1, \dots T_{m+1}$. Each permutation $\sigch{i}$ of the chain corresponding to $T$ is given by a specific traversal of $T$:
\begin{itemize}
\item traverse $T_1, \dots T_{i}$,
\item traverse the root,
\item traverse $T_{i+1}, \dots T_{m+1}$.
\end{itemize}
By doing so, we indeed obtain the same tree-inversions. Beside, if $c = (\sigch{m}, \dots, \sigch{1})$ is the chain corresponding to the maximal class element $\sigma$, then $\sigch{i}$ is the subword formed by the $i^{th}$ occurrences of the letters of $\sigma$.  An example is given in Figure \ref{fig:decreasing-tree-bijection}.

\subsection{Link with $m$-Tamari lattices}
\label{sub-sec:tam}

In \cite{PBT2}, the authors define the \emph{sylvester} congruence on words. As we said in Section \ref{sub-sec:lattice-def}, the metasylvester congruence is finer than the sylvester one. Indeed, the sylvester congruence is the transitive closure of the relation

\begin{align}
ac \dots b &\equiv ca \dots b, & a \leq b < c.
\end{align}
This means that if two $m$-permuations are congruent for the metasylvester relations, they are also congruent for the sylvester one. In the right weak order on permutations, the sylvester classes form intervals. The maximal elements are the permutations avoiding the pattern $132$. It is shown in particular in \cite{PBT2} that the subposet induced by these maximal elements is the Tamari lattice. And the Tamari lattice is also a quotient lattice of the right weak order by the sylvester relation. 

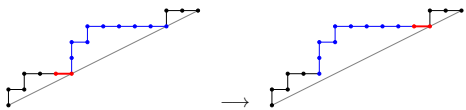
\begin{figure}[ht]
\begin{center}
\scalebox{0.7}{

\begin{tabular}{ccc}

\scalebox{0.3}{
\begin{tikzpicture}
\draw[Gray] (0,0) -- (12,6);
\draw[Path] (0.000, 0.000) -- (0.000, 1.000);
\draw[Path] (0.000, 1.000) -- (1.000, 1.000);
\draw[Path] (1.000, 1.000) -- (1.000, 2.000);
\draw[Path] (1.000, 2.000) -- (2.000, 2.000);
\draw[Path] (2.000, 2.000) -- (3.000, 2.000);
\draw[DPoint] (0.000, 0.000) circle;
\draw[DPoint] (0.000, 1.000) circle;
\draw[DPoint] (1.000, 1.000) circle;
\draw[DPoint] (1.000, 2.000) circle;
\draw[DPoint] (2.000, 2.000) circle;
\draw[StrongPath,Red] (3.000, 2.000) -- (4.000, 2.000);
\draw[DPoint,Red] (3.000, 2.000) circle;
\draw[Path,Blue] (4.000, 2.000) -- (4.000, 3.000);
\draw[DPoint,Red] (4.000, 2.000) circle;
\draw[Path,Blue] (4.000, 3.000) -- (4.000, 4.000);
\draw[Path,Blue] (4.000, 4.000) -- (5.000, 4.000);
\draw[Path,Blue] (5.000, 4.000) -- (5.000, 5.000);
\draw[Path,Blue] (5.000, 5.000) -- (6.000, 5.000);
\draw[Path,Blue] (6.000, 5.000) -- (7.000, 5.000);
\draw[Path,Blue] (7.000, 5.000) -- (8.000, 5.000);
\draw[Path,Blue] (8.000, 5.000) -- (9.000, 5.000);
\draw[Path,Blue] (9.000, 5.000) -- (10.000, 5.000);
\draw[DPoint,Blue] (4.000, 3.000) circle;
\draw[DPoint,Blue] (4.000, 4.000) circle;
\draw[DPoint,Blue] (5.000, 4.000) circle;
\draw[DPoint,Blue] (5.000, 5.000) circle;
\draw[DPoint,Blue] (6.000, 5.000) circle;
\draw[DPoint,Blue] (7.000, 5.000) circle;
\draw[DPoint,Blue] (8.000, 5.000) circle;
\draw[DPoint,Blue] (9.000, 5.000) circle;
\draw[Path] (10.000, 5.000) -- (10.000, 6.000);
\draw[DPoint,Blue] (10.000, 5.000) circle;
\draw[Path] (10.000, 6.000) -- (11.000, 6.000);
\draw[Path] (11.000, 6.000) -- (12.000, 6.000);
\draw[DPoint] (10.000, 6.000) circle;
\draw[DPoint] (11.000, 6.000) circle;
\draw[DPoint] (12.000, 6.000) circle;
\end{tikzpicture}
} &
$\longrightarrow$
&
\scalebox{0.3}{
\begin{tikzpicture}
\draw[Gray] (0,0) -- (12,6);
\draw[Path] (0.000, 0.000) -- (0.000, 1.000);
\draw[Path] (0.000, 1.000) -- (1.000, 1.000);
\draw[Path] (1.000, 1.000) -- (1.000, 2.000);
\draw[Path] (1.000, 2.000) -- (2.000, 2.000);
\draw[Path] (2.000, 2.000) -- (3.000, 2.000);
\draw[DPoint] (0.000, 0.000) circle;
\draw[DPoint] (0.000, 1.000) circle;
\draw[DPoint] (1.000, 1.000) circle;
\draw[DPoint] (1.000, 2.000) circle;
\draw[DPoint] (2.000, 2.000) circle;
\draw[Path,Blue] (3.000, 2.000) -- (3.000, 3.000);
\draw[Path,Blue] (3.000, 3.000) -- (3.000, 4.000);
\draw[Path,Blue] (3.000, 4.000) -- (4.000, 4.000);
\draw[Path,Blue] (4.000, 4.000) -- (4.000, 5.000);
\draw[Path,Blue] (4.000, 5.000) -- (5.000, 5.000);
\draw[Path,Blue] (5.000, 5.000) -- (6.000, 5.000);
\draw[Path,Blue] (6.000, 5.000) -- (7.000, 5.000);
\draw[Path,Blue] (7.000, 5.000) -- (8.000, 5.000);
\draw[Path,Blue] (8.000, 5.000) -- (9.000, 5.000);
\draw[DPoint,Blue] (3.000, 2.000) circle;
\draw[DPoint,Blue] (3.000, 3.000) circle;
\draw[DPoint,Blue] (3.000, 4.000) circle;
\draw[DPoint,Blue] (4.000, 4.000) circle;
\draw[DPoint,Blue] (4.000, 5.000) circle;
\draw[DPoint,Blue] (5.000, 5.000) circle;
\draw[DPoint,Blue] (6.000, 5.000) circle;
\draw[DPoint,Blue] (7.000, 5.000) circle;
\draw[DPoint,Blue] (8.000, 5.000) circle;
\draw[StrongPath,Red] (9.000, 5.000) -- (10.000, 5.000);
\draw[DPoint,Red] (9.000, 5.000) circle;
\draw[Path] (10.000, 5.000) -- (10.000, 6.000);
\draw[DPoint,Red] (10.000, 5.000) circle;
\draw[Path] (10.000, 6.000) -- (11.000, 6.000);
\draw[Path] (11.000, 6.000) -- (12.000, 6.000);
\draw[DPoint] (10.000, 6.000) circle;
\draw[DPoint] (11.000, 6.000) circle;
\draw[DPoint] (12.000, 6.000) circle;
\end{tikzpicture}
}
\end{tabular}
}
\end{center}
\caption{Tamari Rotation on $m$-ballot paths}
\label{fig:mpath-rotation}
\end{figure}

\begin{figure}[ht]
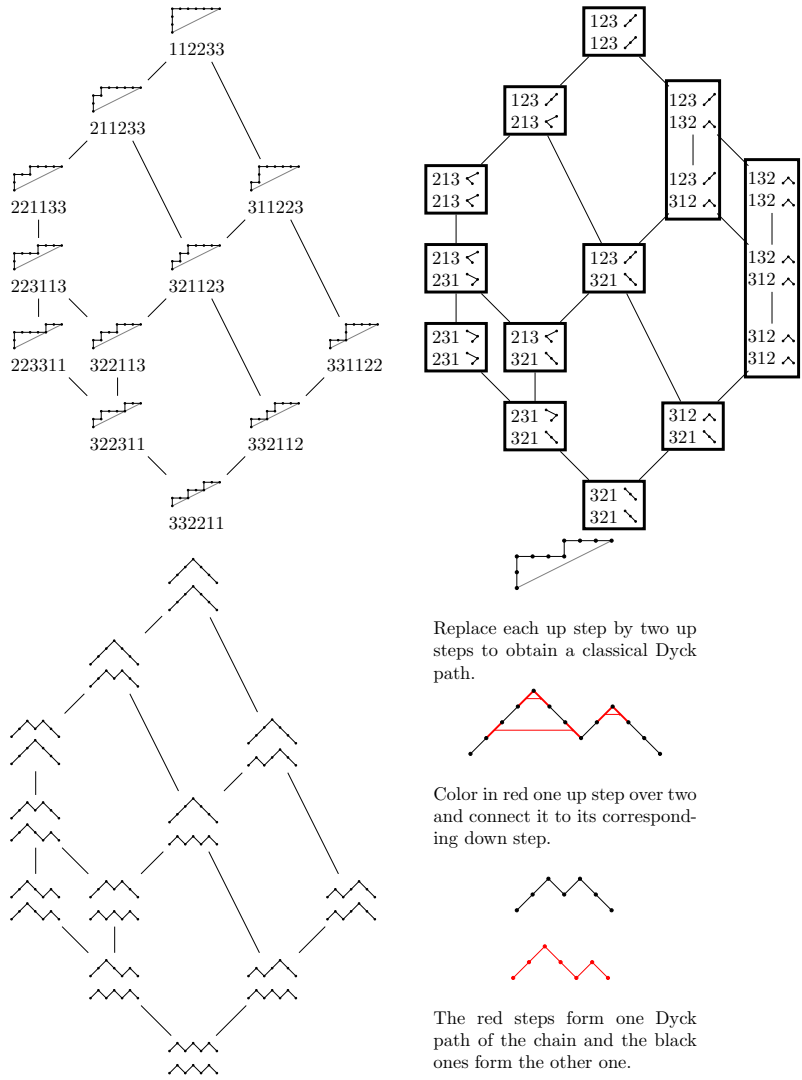

\scalebox{0.7}{
\begin{tabular}{cc}
\begin{minipage}{3in}

\def \wlev{1.5}
\def \hlev{1.5}
\def \sscale{0.15}
\def \fpath{figures/mpaths/}

\begin{tikzpicture}[every text node part/.style={align=center}]
\node(P112233) at (0,0){
\scalebox{\sscale}{\input{\fpath P3-2-l}}\\
112233};
\node(P211233) at (-1 * \wlev, -1 * \hlev){
\scalebox{\sscale}{\input{\fpath P3-2-k}}\\
211233};

\node(P221133) at (-2 * \wlev, -2 * \hlev){
\scalebox{\sscale}{\input{\fpath P3-2-i}}\\
221133};
\node(P311223) at ( 1 * \wlev, -2 * \hlev){
\scalebox{\sscale}{\input{\fpath P3-2-j}}\\
311223};

\node(P223113) at (-2 * \wlev, -3 * \hlev){
\scalebox{\sscale}{\input{\fpath P3-2-g}}\\
223113};
\node(P321123) at ( 0 * \wlev, -3 * \hlev){
\scalebox{\sscale}{\input{\fpath P3-2-h}}\\
321123};

\node(P223311) at (-2 * \wlev, -4 * \hlev){
\scalebox{\sscale}{\input{\fpath P3-2-d}}\\
223311};
\node(P322113) at (-1 * \wlev, -4 * \hlev){
\scalebox{\sscale}{\input{\fpath P3-2-e}}\\
322113};
\node(P331122) at ( 2 * \wlev, -4 * \hlev){
\scalebox{\sscale}{\input{\fpath P3-2-f}}\\
331122};

\node(P322311) at (-1 * \wlev, -5 * \hlev){
\scalebox{\sscale}{\input{\fpath P3-2-b}}\\
322311};
\node(P332112) at ( 1 * \wlev, -5 * \hlev){
\scalebox{\sscale}{\input{\fpath P3-2-c}}\\
332112};

\node(P332211) at ( 0 * \wlev, -6 * \hlev){
\scalebox{\sscale}{\input{\fpath P3-2-a}}\\
332211};

\draw (P112233) -- (P211233);
\draw (P112233) -- (P311223);

\draw(P211233) -- (P221133);
\draw(P211233) -- (P321123);

\draw(P221133) -- (P223113);
\draw(P311223) -- (P321123);
\draw(P311223) -- (P331122);

\draw(P223113) -- (P223311);
\draw(P223113) -- (P322113);
\draw(P321123) -- (P322113);
\draw(P321123) -- (P332112);

\draw(P223311) -- (P322311);
\draw(P322113) -- (P322311);
\draw(P331122) -- (P332112);

\draw(P322311) -- (P332211);
\draw(P332112) -- (P332211);
\end{tikzpicture}
\end{minipage} &
\begin{minipage}{3in}
\def \wlev{1.5}
\def \hlev{1.5}
\def \sscale{0.2}
\def \fpath{figures/binary/}

\tikzstyle{TamNode} = [rectangle, draw, ultra thick]

\begin{tikzpicture}[every text node part/.style={align=center}]

\node[TamNode](P123_123) at (0,0){
123 \scalebox{\sscale}{\input{\fpath T3-1}} \\
123 \scalebox{\sscale}{\input{\fpath T3-1}}
};

\node[TamNode](P123_213) at (-1 * \wlev, -1 * \hlev){
123 \scalebox{\sscale}{\input{\fpath T3-1}} \\
213 \scalebox{\sscale}{\input{\fpath T3-2}}
};

\node[TamNode, fit = (P123_132) (P123_312)](fit1){};
\node(P123_132) at (1 * \wlev, -1 * \hlev){
123 \scalebox{\sscale}{\input{\fpath T3-1}} \\
132 \scalebox{\sscale}{\input{\fpath T3-4}}
};

\node[TamNode](P213_213) at (-2 * \wlev, -2 * \hlev){
213 \scalebox{\sscale}{\input{\fpath T3-2}} \\
213 \scalebox{\sscale}{\input{\fpath T3-2}}
};

\node(P123_312) at (1 * \wlev, -2 * \hlev){
123 \scalebox{\sscale}{\input{\fpath T3-1}} \\
312 \scalebox{\sscale}{\input{\fpath T3-4}}
};

\node[TamNode, fit = (P132_132) (P132_312) (P312_312)](fit2){};
\node(P132_132) at (2 * \wlev, -2 * \hlev){
132 \scalebox{\sscale}{\input{\fpath T3-4}} \\
132 \scalebox{\sscale}{\input{\fpath T3-4}}
};

\node[TamNode](P213_231) at (-2 * \wlev, -3 * \hlev){
213 \scalebox{\sscale}{\input{\fpath T3-2}} \\
231 \scalebox{\sscale}{\input{\fpath T3-3}}
};

\node[TamNode](P123_321) at (0 * \wlev, -3 * \hlev){
123 \scalebox{\sscale}{\input{\fpath T3-1}} \\
321 \scalebox{\sscale}{\input{\fpath T3-5}}
};

\node(P132_312) at (2 * \wlev, -3 * \hlev){
132 \scalebox{\sscale}{\input{\fpath T3-4}} \\
312 \scalebox{\sscale}{\input{\fpath T3-4}}
};

\node[TamNode](P231_231) at (-2 * \wlev, -4 * \hlev){
231 \scalebox{\sscale}{\input{\fpath T3-3}} \\
231 \scalebox{\sscale}{\input{\fpath T3-3}}
};

\node[TamNode](P213_321) at (-1 * \wlev, -4 * \hlev){
213 \scalebox{\sscale}{\input{\fpath T3-2}} \\
321 \scalebox{\sscale}{\input{\fpath T3-5}}
};

\node(P312_312) at (2 * \wlev, -4 * \hlev){
312 \scalebox{\sscale}{\input{\fpath T3-4}} \\
312 \scalebox{\sscale}{\input{\fpath T3-4}}
};

\node[TamNode](P231_321) at (-1 * \wlev, -5 * \hlev){
231 \scalebox{\sscale}{\input{\fpath T3-3}} \\
321 \scalebox{\sscale}{\input{\fpath T3-5}}
};

\node[TamNode](P312_321) at (1 * \wlev, -5 * \hlev){
312 \scalebox{\sscale}{\input{\fpath T3-4}} \\
321 \scalebox{\sscale}{\input{\fpath T3-5}}
};

\node[TamNode](P321_321) at (0 * \wlev, -6 * \hlev){
321 \scalebox{\sscale}{\input{\fpath T3-5}} \\
321 \scalebox{\sscale}{\input{\fpath T3-5}}
};

\draw (P123_123) -- (P123_213);
\draw (P123_123) -- (P123_132);
\draw (P123_213) -- (P213_213);
\draw (P123_213) -- (P123_321);
\draw (P123_132) -- (P123_312);
\draw (P123_132) -- (P132_132);
\draw (P213_213) -- (P213_231);
\draw (P123_312) -- (P123_321);
\draw (P123_312) -- (P132_312);
\draw (P132_132) -- (P132_312);
\draw (P213_231) -- (P231_231);
\draw (P213_231) -- (P213_321);
\draw (P123_321) -- (P213_321);
\draw (P123_321) -- (P312_321);
\draw (P132_312) -- (P312_312);
\draw (P231_231) -- (P231_321);
\draw (P213_321) -- (P231_321);
\draw (P312_312) -- (P312_321);
\draw (P231_321) -- (P321_321);
\draw (P312_321) -- (P321_321);

\end{tikzpicture}
\end{minipage}
\\
\begin{minipage}{3in}

\def \wlev{1.5}
\def \hlev{1.5}
\def \sscale{0.15}
\def \fpath{figures/dyck/}

\begin{tikzpicture}[every text node part/.style={align=center}]
\node(P112233) at (0,0){
\scalebox{\sscale}{\input{\fpath D3-5}}\\
\scalebox{\sscale}{\input{\fpath D3-5}}};

\node(P211233) at (-1 * \wlev, -1 * \hlev){
\scalebox{\sscale}{\input{\fpath D3-5}}\\
\scalebox{\sscale}{\input{\fpath D3-3}}};

\node(P221133) at (-2 * \wlev, -2 * \hlev){
\scalebox{\sscale}{\input{\fpath D3-3}}\\
\scalebox{\sscale}{\input{\fpath D3-5}}};

\node(P311223) at ( 1 * \wlev, -2 * \hlev){
\scalebox{\sscale}{\input{\fpath D3-5}}\\
\scalebox{\sscale}{\input{\fpath D3-4}}};

\node(P223113) at (-2 * \wlev, -3 * \hlev){
\scalebox{\sscale}{\input{\fpath D3-3}}\\
\scalebox{\sscale}{\input{\fpath D3-2}}};

\node(P321123) at ( 0 * \wlev, -3 * \hlev){
\scalebox{\sscale}{\input{\fpath D3-5}}\\
\scalebox{\sscale}{\input{\fpath D3-1}}};

\node(P223311) at (-2 * \wlev, -4 * \hlev){
\scalebox{\sscale}{\input{\fpath D3-2}}\\
\scalebox{\sscale}{\input{\fpath D3-2}}};

\node(P322113) at (-1 * \wlev, -4 * \hlev){
\scalebox{\sscale}{\input{\fpath D3-3}}\\
\scalebox{\sscale}{\input{\fpath D3-1}}};

\node(P331122) at ( 2 * \wlev, -4 * \hlev){
\scalebox{\sscale}{\input{\fpath D3-4}}\\
\scalebox{\sscale}{\input{\fpath D3-4}}};

\node(P322311) at (-1 * \wlev, -5 * \hlev){
\scalebox{\sscale}{\input{\fpath D3-2}}\\
\scalebox{\sscale}{\input{\fpath D3-1}}};

\node(P332112) at ( 1 * \wlev, -5 * \hlev){
\scalebox{\sscale}{\input{\fpath D3-4}}\\
\scalebox{\sscale}{\input{\fpath D3-1}}};

\node(P332211) at ( 0 * \wlev, -6 * \hlev){
\scalebox{\sscale}{\input{\fpath D3-1}}\\
\scalebox{\sscale}{\input{\fpath D3-1}}};

\draw (P112233) -- (P211233);
\draw (P112233) -- (P311223);

\draw(P211233) -- (P221133);
\draw(P211233) -- (P321123);

\draw(P221133) -- (P223113);
\draw(P311223) -- (P321123);
\draw(P311223) -- (P331122);

\draw(P223113) -- (P223311);
\draw(P223113) -- (P322113);
\draw(P321123) -- (P322113);
\draw(P321123) -- (P332112);

\draw(P223311) -- (P322311);
\draw(P322113) -- (P322311);
\draw(P331122) -- (P332112);

\draw(P322311) -- (P332211);
\draw(P332112) -- (P332211);
\end{tikzpicture}
\end{minipage}
&
\begin{minipage}{3in}
\begin{tabular}{m{5cm}}
\begin{center}
\scalebox{0.3}{

\begin{tikzpicture}
\draw[Gray] (0,0) -- (6,3);
\draw (0.000, 0.000) -- (0.000, 1.000);
\draw (0.000, 1.000) -- (0.000, 2.000);
\draw (0.000, 2.000) -- (1.000, 2.000);
\draw (1.000, 2.000) -- (2.000, 2.000);
\draw (2.000, 2.000) -- (3.000, 2.000);
\draw (3.000, 2.000) -- (3.000, 3.000);
\draw (3.000, 3.000) -- (4.000, 3.000);
\draw (4.000, 3.000) -- (5.000, 3.000);
\draw (5.000, 3.000) -- (6.000, 3.000);
\draw[DPoint] (0.000, 0.000) circle;
\draw[DPoint] (0.000, 1.000) circle;
\draw[DPoint] (0.000, 2.000) circle;
\draw[DPoint] (1.000, 2.000) circle;
\draw[DPoint] (2.000, 2.000) circle;
\draw[DPoint] (3.000, 2.000) circle;
\draw[DPoint] (3.000, 3.000) circle;
\draw[DPoint] (4.000, 3.000) circle;
\draw[DPoint] (5.000, 3.000) circle;
\draw[DPoint] (6.000, 3.000) circle;
\end{tikzpicture}}
\end{center}
 \\

Replace each up step by two up steps to obtain a classical Dyck path.
\\
\begin{center}
\scalebox{0.3}{
\begin{tikzpicture}
\draw[Path] (0.000, 0.000) -- (1.000, 1.000);
\draw[RedPath, StrongPath] (1.000, 1.000) -- (2.000, 2.000);
\draw[Path] (2.000, 2.000) -- (3.000, 3.000);
\draw[RedPath, StrongPath] (3.000, 3.000) -- (4.000, 4.000);
\draw[RedPath, StrongPath] (4.000, 4.000) -- (5.000, 3.000);
\draw[Path] (5.000, 3.000) -- (6.000, 2.000);
\draw[RedPath, StrongPath] (6.000, 2.000) -- (7.000, 1.000);
\draw[Path] (7.000, 1.000) -- (8.000, 2.000);
\draw[RedPath, StrongPath] (8.000, 2.000) -- (9.000, 3.000);
\draw[RedPath, StrongPath] (9.000, 3.000) -- (10.000, 2.000);
\draw[Path] (10.000, 2.000) -- (11.000, 1.000);
\draw[Path] (11.000, 1.000) -- (12.000, 0.000);
\draw[DPoint] (0.000, 0.000) circle;
\draw[DPoint] (1.000, 1.000) circle;
\draw[DPoint] (2.000, 2.000) circle;
\draw[DPoint] (3.000, 3.000) circle;
\draw[DPoint] (4.000, 4.000) circle;
\draw[DPoint] (5.000, 3.000) circle;
\draw[DPoint] (6.000, 2.000) circle;
\draw[DPoint] (7.000, 1.000) circle;
\draw[DPoint] (8.000, 2.000) circle;
\draw[DPoint] (9.000, 3.000) circle;
\draw[DPoint] (10.000, 2.000) circle;
\draw[DPoint] (11.000, 1.000) circle;
\draw[DPoint] (12.000, 0.000) circle;

\draw[RedPath, -> ] (1.5,1.5) -- (6.5,1.5);
\draw[RedPath, -> ] (3.5,3.5) -- (4.5,3.5);
\draw[RedPath, -> ] (8.5,2.5) -- (9.5,2.5);
\end{tikzpicture}
}
\end{center}

\\

Color in red one up step over two and connect it to its corresponding down step.

\\
\begin{center}
\scalebox{0.3}{


\begin{tikzpicture}
\draw[Path] (0.000, 0.000) -- (1.000, 1.000);
\draw[Path] (1.000, 1.000) -- (2.000, 2.000);
\draw[Path] (2.000, 2.000) -- (3.000, 1.000);
\draw[Path] (3.000, 1.000) -- (4.000, 2.000);
\draw[Path] (4.000, 2.000) -- (5.000, 1.000);
\draw[Path] (5.000, 1.000) -- (6.000, 0.000);
\draw[DPoint] (0.000, 0.000) circle;
\draw[DPoint] (1.000, 1.000) circle;
\draw[DPoint] (2.000, 2.000) circle;
\draw[DPoint] (3.000, 1.000) circle;
\draw[DPoint] (4.000, 2.000) circle;
\draw[DPoint] (5.000, 1.000) circle;
\draw[DPoint] (6.000, 0.000) circle;
\end{tikzpicture}}
\end{center}
\\
\begin{center}
\red{
\scalebox{0.3}{


\begin{tikzpicture}
\draw[Path] (0.000, 0.000) -- (1.000, 1.000);
\draw[Path] (1.000, 1.000) -- (2.000, 2.000);
\draw[Path] (2.000, 2.000) -- (3.000, 1.000);
\draw[Path] (3.000, 1.000) -- (4.000, 0.000);
\draw[Path] (4.000, 0.000) -- (5.000, 1.000);
\draw[Path] (5.000, 1.000) -- (6.000, 0.000);
\draw[DPoint] (0.000, 0.000) circle;
\draw[DPoint] (1.000, 1.000) circle;
\draw[DPoint] (2.000, 2.000) circle;
\draw[DPoint] (3.000, 1.000) circle;
\draw[DPoint] (4.000, 0.000) circle;
\draw[DPoint] (5.000, 1.000) circle;
\draw[DPoint] (6.000, 0.000) circle;
\end{tikzpicture}}
}
\end{center}
\\

The red steps form one Dyck path of the chain and the black ones form the other one.

\end{tabular}
\end{minipage}
\end{tabular}
}
\caption{On the top: the $m$-Tamari lattice as a sublattice and quotient lattice of the metasylvester lattice. On the bottom: the realization of $m$-Tamari in terms of chains of Dyck paths.}
\label{fig:mtam}
\end{figure}

The $m$-Tamari lattice is a generalization of the Tamari lattice which was introduced in \cite{BergmTamari}. It is traditionally defined on $m$-ballot paths : a path from $(0,0)$ to $(nm,n)$ made from horizontal steps $(1,0)$ and vertical steps $(0,1)$ which always stays above the line $y=\frac{x}{m}$. The cover relation is given by the rotation operation illustrated in Figure \ref{fig:mpath-rotation}. In \cite{mTamari}, it was shown that the $m$-Tamari lattice of size $n$, denoted by $\Tamnm$, can be understood as an upper ideal of the classical Tamari lattice of size $m \times n$. For a detailed definition of the $m$-Tamari lattice and its properties, we invite our lectors to read \cite{mTamari}. The relations between $m$-permutations, $m$-Tamari lattices and trees is explained in \cite{mTamAlg} following some combinatorial constructions given in \cite{ME_THESE}. In this paper, we are interested in the definition of the $m$-Tamari lattice in terms of the sylvester relation on $m$-permutations. As explained in Section \ref{sub-sec:mperms}, the lattice on $m$-permutations of size $n$ is a lower ideal of the right weak order of size $n \times m$. It corresponds to a lower ideal of the Tamari lattice of size $n \times m$, and more precisely, it is the symmetric of the $m$-Tamari ideal. In other words, the $m$-Tamari lattice is both a sublattice and a quotient lattice of the lattice on $m$-permutations. More precisely, if $\sigma$ and $\mu$ are two $m$-permtuations which are maximal elements of sylvester classes, then $\sigma \leq \mu$ in the $m$-permutations lattice if and only if $T_\sigma \geq T_\mu$ in the $m$-Tamari lattice where $T_\sigma$ and $T_\mu$ are the corresponding elements of respectively $\sigma$ and $\mu$, as illustrated in Figure \ref{fig:mtam}.

The following result is direct from the previous remarks.

\begin{Theoreme}
The $m$-Tamari lattice is both a sublattice and a quotient lattice of the metasylvester lattice.
\end{Theoreme}

The explicit surjection from $m$-permutations to $m$-ballot paths use the notion of $m$-binary trees defined in \cite{ME_THESE} and the \emph{binary search tree} insertion algorithm already used in the case of the classical Tamari lattice \cite{PBT2}. We will not detail it here. The realization of the metasylvester lattice on chains of permutations gives us an interesting result on the $m$-Tamari lattice. Indeed, by applying the classical sylvester surjection by binary search tree insertion on all permutations of a metasylvester $m$-chain, we get a Tamari chain of binary trees. These chains give a new realization of the $m$-Tamari lattice. The third lattice of Figure~\ref{fig:mtam} represents the $m$-Tamari lattice in terms of this new realization: we have constructed the binary tree associated to each permutation of the metasylvester chain (as in the second lattice of Figure \ref{fig:mtam}) and then taken its symmetric Dyck path in the Tamari lattice. On the right, we give a succinct description of the direct algorithm from a $m$-ballot path to a chain of Dyck paths. The chains we obtain can be characterized through the same pattern avoidance characterizing metasylvester $m$-chains. Furthermore, there have been some other descriptions of the $m$-Tamari lattice in terms of chains. In particular, our realization seems to correspond to the lattice closure of the $m$-cover poset of the Tamari lattice described in \cite{mCover}.

\begin{footnotesize}
\bibliographystyle{plain}
\bibliography{../../../biblio/all}
\end{footnotesize}

\end{document}